\documentclass[amstex,12pt,russian,amssymb]{article}

\usepackage{mathtext}
\usepackage[cp1251]{inputenc}
\usepackage[T2A]{fontenc}
\usepackage[russian]{babel}
\usepackage[dvips]{graphicx}
\usepackage{amsmath}
\usepackage{amssymb}
\usepackage{amsxtra}
\usepackage{latexsym}
\usepackage{ifthen}

\begin{document}

\newcounter{lemma}
\newcommand{\lemma}{\par \refstepcounter{lemma}%
{\bf Лемма \arabic{lemma}.}}

\newcounter{corollary}
\newcommand{\corollary}{\par \refstepcounter{corollary}%
{\bf Следствие \arabic{corollary}.}}

\newcounter{remark}
\newcommand{\remark}{\par \refstepcounter{remark}%
{\bf Замечание \arabic{remark}.}}

\newcounter{theorem}
\newcommand{\theorem}{\par \refstepcounter{theorem}%
{\bf Теорема \arabic{theorem}.}}

\newcounter{proposition}
\newcommand{\proposition}{\par \refstepcounter{proposition}%
{\bf Предложение \arabic{proposition}.}}

\renewcommand{\refname}{\centerline{\bf Список литературы}}

\newcommand{\proof}{{\it Доказательство.\,\,}}

\noindent УДК 517.5

{\bf Е.А.~Севостьянов}

{\bf Є.О.~Севостьянов}

{\bf E.A.~Sevost'yanov}

{\small Институт прикладной математики и механики НАН Украины,
Донецк}

{\small Інститут прикладної математики і механіки НАН України,
Донецьк}

{\small Institute of Applied Mathematics and Mechanics of NAS of
Ukraine, Donetsk}

{\bf Аналог теоремы Монтеля для отображений класса Соболева с
конечным искажением}

{\bf Аналог теореми Монтеля для відображень класу Соболєва зі
скінченним спотворенням}

{\bf Analog of Montel theorem for mappings of Sobolev class with
finite distortion}

\medskip
Настоящая работа посвящена изучению классов отображений с
неограниченной характеристикой квазиконформности. Получен результат
о нормальности семейств открытых дискретных отображений
$f:D\rightarrow {\Bbb C}\setminus\{a, b\}$ класса $W_{loc}^{1, 1},$
имеющих конечное искажение и не принимающих два фиксированных
значения $a\ne b$ в ${\Bbb C},$ максимальная дилатация которых имеет
мажоранту конечного среднего колебания в каждой точке. Указанный
результат справедлив, в частности, для так называемых
$Q$-отображений и является аналогом известной теоремы Монтеля для
аналитических функций.

\medskip
Дану роботу присвячено вивченню класів відображень з необмеженою
характеристикою квазіконформності. Отримано результат про
нормальність сімей відкритих дискретних відображень $f:D\rightarrow
{\Bbb C}\setminus\{a, b\}$ класу $W_{loc}^{1, 1},$ що мають
скінченне спотворення і не приймають принаймні два фіксованих
значення $a\ne b$ в ${\Bbb C},$ максимальна дилатація котрих має
мажоранту скінченного середнього коливання в кожній точці. Вказаний
результат є справедливим, зокрема, для так званих $Q$-відображень і
є аналогом відомої теореми Монтеля для аналітичних функцій.

\medskip
The present paper is devoted to the study of classes of mappings
with non-bounded characteristic of quasiconformality. It is obtained
a result on normal families of the open discrete mappings
$f:D\rightarrow {\Bbb C}\setminus\{a, b\}$ of the class $W_{loc}^{1,
1}$ having a finite distortion and omitting two fixed values $a\ne
b$ in ${\Bbb C},$ maximal dilatations of which has a majorant of the
class of finite mean oscillation at every point. In particular, the
result mentioned above holds for the so-called $Q$-mappings and is
an analog of known Montel theorem for analytic functions.

\newpage

{\bf 1. Введение.} Настоящая заметка посвящена обобщению одного
аналога известной теоремы Монтеля о нормальности семейств
аналитических функций (см. \cite[$\S$ 32, гл. II]{Mont}). Ввиду этой
теоремы, как известно, семейство $\frak{F}_{a, b}(D)$ аналитических
функций $f:D\rightarrow {\Bbb C}\setminus\{a, b\}$ области $D\subset
{\Bbb C}$ является нормальным при всяких фиксированных значениях $a,
b\in {\Bbb C},$ $a\ne b$ (см. там же). Как оказалось, указанный
результат остаётся также справедливым и для более общих классов
открытых дискретных отображений класса Соболева $W_{loc}^{1,1}$ с
конечным искажением, как только так называемая дилатация
$K_{\mu}(z)$ этих отображений удовлетворяет некоторым (достаточно
общим) ограничениям на рост.

\medskip
Более того, следует отметить, что указанный результат справедлив
также для так называемых $Q$-отображений, исследованных автором
ранее (см., напр., \cite[раздел~5]{Sev$_1$}). Одно из
сформулированных в данной работе утверждений по этому поводу
усиливает более ранние результаты о нормальности семейств
$Q$-отображений, не принимающих значения множества $E$ положительной
конформной ёмкости. Вместо этого в настоящей заметке предлагается
ограничиться лишь двуточечным множеством $E$ комплексной плоскости.
Следует также отметить, что здесь речь идёт лишь о случае
размерности пространства ${\Bbb R}^n$ при $n=2,$ поскольку случаи
больших размерностей требуют иных подходов нежели те, что
рассмотрены ниже.

\medskip
Приведём необходимые определения и обозначения. Всюду далее $m$ --
мера Лебега в ${\Bbb C},$ $D$ -- область в ${\Bbb C},$
$\overline{\Bbb C}={\Bbb C}\cup\{\infty\}$ -- одноточечная
компактификация ${\Bbb C}.$

\medskip
Для открытого множества $U\subset {\Bbb C}$ символом $C_1^{\,0}(U),$
как обычно, обозначается множество всех непрерывно дифференцируемых
в $U$ функций с компактным носителем в $U.$ Пусть $u:U\rightarrow
{\Bbb R}$ -- некоторая функция, $u\in L_{loc}^{\,1}(U),$
$z=x_1+ix_2,$ $z\in {\Bbb C}.$ Предположим, что найдётся функция
$v\in L_{loc}^{\,1}(U),$ такая что
$$\int\limits_U \frac{\partial \varphi}{\partial x_i}(z)u(z)dm(z)=
-\int\limits_U \varphi(z)v(z)dm(z)$$
для любой функции $\varphi\in C_1^{\,0}(U).$ Тогда будем говорить,
что функция $v$ является {\it обобщённой производной первого порядка
функции $u$ по переменной $x_i$} и обозначать символом:
$\frac{\partial u}{\partial x_i}(z):=v.$ Функция $u\in
W_{loc}^{1,1}(U),$ если $u$ имеет обобщённые производные первого
порядка по каждой из переменных в $U,$ которые являются локально
интегрируемыми в $U.$

\medskip
Пусть $G$ -- открытое множество в ${\Bbb C}.$ Отображение
$f:G\rightarrow {\Bbb C}$ принадлежит {\it классу Соболева}
$W^{1,1}_{loc}(G),$ пишут $f\in W^{1,1}_{loc}(G),$ если его
координатные функции $f_k,$ $k=1, 2,$ $f(z)=f_1(z)+ if_2(z),$
обладают обобщёнными частными производными первого порядка, которые
локально интегрируемы в $G$ в первой степени. Если $f\in
W^{1,1}_{loc}(G),$ то {\it матрицу Якоби} $f^{\,\prime}(z)$
отображения $f$ в точке $z$ определим как матрицу $(2\times 2),$
составленную из формальных частных производных по Соболеву (в $i$-й
строке и $j$-м столбце располагается элемент $\frac{\partial
f_i}{\partial x_j}(z)$). Стоит отметить, что в почти всех точках
дифференцируемости $z\in G$ элементы $\frac{\partial f_i}{\partial
x_j}(z),$ понимаемые в указанном выше смысле, совпадают с обычной
частной производной функции $f_i$ по переменной $x_j,$ также
обозначаемой символом $\frac{\partial f_i}{\partial x_j}(z),$ если
невозможно недоразумение (см. \cite[теорема~1, п.~1.1.3, $\S$~1.1,
гл.~I]{Ma}). Отметим также, что всякое отображение $f\in
W^{1,1}_{loc}(G)$ почти всюду имеет обычные частные производные по
каждой из переменных (см. там же).

\medskip
Для комплекснозначной функции $f:D\rightarrow {\Bbb C},$ заданной в
области $D\subset {\Bbb C},$ имеющей частные производные по $x$ и
$y$ при почти всех $z=x + iy,$ полагаем $\overline{\partial} f=
f_{\overline{z}} = \left(f_x + if_y\right)/2$ и $\partial f = f_z =
\left(f_x - if_y\right)/2.$ Полагаем
$\mu(z)=\mu_f(z)=f_{\overline{z}}/f_z,$ при $f_z \ne 0$ и $\mu(z)=0$
в противном случае. Указанная комплекснозначная функция $\mu$
называется {\it комплексной дилатацией} отображения $f$ в точке $z.$
{\it Мак\-си\-маль\-ной дилатацией} отображения $f$ в точке $z$
называется следующая функция:
\begin{equation}\label{eq1.22A}
K_{\mu_f}(z)\quad=\quad K_{\mu}(z)\quad=\quad\frac{1 + |\mu (z)|}{|1
- |\mu\,(z)||}\,.
\end{equation}

\medskip
Заметим, что $J(f, z)=|f_z|^2-|f_{\overline{z}}|^2,$ где $J(f,
z):={\rm det\,}f^{\,\prime}(z),$
что может быть проверено прямым подсчётом (см., напр.,
\cite[пункт~C, гл.~I]{A}).

\medskip
Отображение $f:D\rightarrow {\Bbb C}$ называется {\it отображением с
конечным искажением}, если $f\in W_{loc}^{1,1}(D)$ и для некоторой
функции $K(z): D\rightarrow [1,\infty)$ выполнено условие
$$\Vert f^{\,\prime}(z)\Vert^{2}\le K(z)\cdot
|J(f, z)|$$
при почти всех $z\in D,$ где $\Vert
f^{\,\prime}(z)\Vert=|f_z|+|f_{\overline{z}}|$ (см. \cite[п.~6.3,
гл.~VI]{IM}). Суть понятия отображения с конечным искажением
заключается в том, что у указанного отображения $f$ матричная норма
производной $\Vert f^{\,\prime}(z)\Vert$ равна нулю в почти всех
точках вырождения якобиана $J(f, z).$

\medskip
Пусть $\left(X,\,d\right)$ и
$\left(X^{\,{\prime}},{d}^{\,{\prime}}\right)$ -- метрические
пространства с расстояниями  $d$  и ${d}^{\,{\prime}}$
соответственно. Семейство $\frak{F}$ непрерывных отображений
$f:X\rightarrow {X}^{\,\,\prime}$ называется {\it нормальным}, если
из любой последовательности отображений $f_{m} \in \frak{F}$ можно
выделить подпоследовательность $f_{m_{k}}$, которая сходится
локально равномерно в $X$ к непрерывной функции $f:X\rightarrow
X^{\,\prime}.$ Отметим, что всюду далее, если не оговорено
противное, $(X, d)=(D, |\cdot|),$ где $D$ -- область в ${\Bbb C},$ а
$|\cdot|$ -- евклидова метрика, $|x-y|=\sqrt{\sum\limits_{i=1}^2
(y_i-x_i)^2},$ где $x=x_1+ix_2,$ $y=y_1+iy_2);$
$\left(X^{\,\prime},\, d^{\,\prime}\right)=\left(\overline{\Bbb
C},\, h\right),$ где $h$ -- хордальная метрика,
$$h(x,\infty)=\frac{1}{\sqrt{1+{|x|}^2}}\,,\quad
h(x,y)=\frac{|x-y|}{\sqrt{1+{|x|}^2} \sqrt{1+{|y|}^2}}\,,  x\ne
\infty\ne y\,.
$$

\medskip
Определение и примеры функций $\varphi:D\rightarrow {\Bbb R}$ класса
$FMO(z_0),$ $z_0\in {\Bbb C},$ являющихся обобщением класса $BMO,$
см., напр., в работе \cite{IR}. Пусть $Q:D\rightarrow [0,\infty]$ --
измеримая по Лебегу функция, тогда $q_{z_0}(r)$ означает среднее
интегральное значение $Q(x)$ над окружностью $|z-z_0|=r,$
\begin{equation}\label{eq32*}
q_{z_0}(r):=\frac{1}{2\pi r}\int\limits_{|z-z_0|=r}Q(z)\,d{\mathcal
H}^{1}\,,
\end{equation}
где ${\mathcal H}^{1}$ -- $1$-мерная мера Хаусдорфа.

\medskip
Для фиксированных области $D\subset {\Bbb C},$ чисел $a, b\in {\Bbb
C},$ $a\ne b,$ и измеримой по Лебегу функции $Q:D\rightarrow [0,
\infty],$ обозначим символом $\frak{G}_{a, b, Q}(D)$ семейство всех
открытых дискретных отображений $f:D\rightarrow {\Bbb
C}\setminus\{a, b\}$ класса $W_{loc}^{1,1}(D)$ и имеющих конечное
искажение, таких что $K_{\mu_f}(z)\le Q(z)$ при почти всех $z\in D.$
Один из основных результатов настоящей работы может быть
сформулирован следующим образом.

\medskip
\begin{theorem}\label{th1A}
{\sl Семейство отображений $\frak{G}_{a, b, Q}(D)$ является
нормальным, как только выполнено, по крайней мере, одно из следующих
условий:
1) $Q\in FMO(z_0)$ в каждой точке $z_0\in D;$ 2) $q_{z_0}(r)\le
C(z_0)\cdot \log\frac{1}{r}$ при $r\rightarrow 0$ и каждой точке
$z_0\in D,$ где $C(z_0)>0$ -- некоторая постоянная; 3) $Q\in
L_{loc}^1(D)$ и при некотором $\varepsilon_0=\varepsilon_0(z_0)$
имеет место соотношение
\begin{equation}\label{eq12}
\int\limits_{0}^{\varepsilon_0} \frac{dt}{tq_{z_0}(t)}=\infty\,.
\end{equation}
}
\end{theorem}

\medskip
Приведём ещё результат по этому поводу. Далее $M$ обозначает
конформный модуль семейства кривых (см., напр., \cite[разд.~6,
гл.~I]{Va$_1$}). Пусть $Q:D\rightarrow [0,\infty]$ -- некоторая
фиксированная вещественнозначная функция. Согласно
\cite[гл.~4]{MRSY}, отображение $f:D\rightarrow \overline{\Bbb C}$
условимся называть {\it $Q$-отоб\-ра\-же\-ни\-ем}, если $f$
удовлетворяет соотношению
$$M(f(\Gamma))\le \int\limits_D Q(z)\cdot \rho^2 (z) dm(z)$$
для произвольного семейства кривых $\Gamma$ в области $D$ и каждой
допустимой функции $\rho\in {\rm adm\,}\Gamma.$ (Определение условия
допустимости $\rho\in {\rm adm\,}\Gamma$ функции $\rho$ см. в
\cite[гл.~I]{Va$_1$}). В частности, если $f$ -- гомеоморфизм, будем
называть такое отображение $Q$-{\it го\-ме\-о\-мор\-физ\-мом}.

\medskip
Для фиксированных области $D\subset {\Bbb C},$ чисел $a, b\in {\Bbb
C},$ $a\ne b,$ и измеримой по Лебегу функции $B:D\rightarrow [0,
\infty],$ обозначим символом $\frak{G^{\,*}}_{a, b, B}(D)$ семейство
всех открытых дискретных $Q$-отображений $f:D\rightarrow {\Bbb
C}\setminus\{a, b\}$ таких, что $Q(z)\le B(z)$ при почти всех $z\in
D.$ Имеет место следующее

\medskip
\begin{corollary}\label{cor1A}
{\sl Семейство отображений $\frak{G^{\,*}}_{a, b, B}(D)$ является
нормальным, как только выполнено, по крайней мере, одно из следующих
условий:
1) $B\in FMO(z_0)$ в каждой точке $z_0\in D;$ 2) $b_{z_0}(r)\le
C(z_0)\cdot \log\frac{1}{r}$ при $r\rightarrow 0$ и каждой точке
$z_0\in D,$ где $C(z_0)>0$ -- некоторая постоянная, а $b_{z_0}(r)$
-- среднее значение функции $B(z)$ над окружностью $S(z_0, r);$ 3)
$B\in L_{loc}^1(D)$ и при некотором
$\varepsilon_0=\varepsilon_0(z_0)$ имеет место соотношение
(\ref{eq12}), где вместо $q_{z_0}(r)$ следует взять $b_{z_0}(r)$ --
среднее значение функции $B(z)$ над окружностью $S(z_0, r).$}
\end{corollary}

\medskip
{\bf 3. Формулировка и доказательство основной леммы.} Докажем,
прежде всего, следующее утверждение.

\medskip
\begin{lemma}\label{lem1}
{\sl Пусть $f:D\rightarrow {\Bbb C}$ -- отображение с конечным
искажением, $f\in W_{loc}^{1,1},$ имеющее вид $f=\varphi\circ g,$
где $g$ -- некоторый гомеоморфизм, а $\varphi$ -- аналитическая
функция. Тогда также $g\in W_{loc}^{1,1}$ и, кроме того, $g$ имеет
конечное искажение}.
\end{lemma}

\medskip
\begin{proof} Пусть $f=\varphi\circ g,$
где $g$ -- некоторый гомеоморфизм, а $\varphi$ -- аналитическая
функция, при этом, $f\in W_{loc}^{1,1}$ и $f$ имеет конечное
искажение. Отметим, что множество точек ветвления
$B_{\varphi}\subset g(D)$ функции $\varphi$ состоит только из
изолированных точек (см. \cite[пункты 5 и 6 (II), гл.~V]{St}).
Следовательно, $g(z)=\varphi^{-1}\circ f$ локально, вне множества
$g^{-1}\left(B_{\varphi}\right).$ Ясно, что множество
$g^{-1}\left(B_{\varphi}\right)$ также состоит из изолированных
точек, следовательно, $g\in ACL(D)$ как композиция аналитической
функции $\varphi^{-1}$ и отображения $f\in W_{loc}^{1,1}(D).$

Покажем, что $g\in W_{loc}^{1,1}(D).$ Пусть далее $\mu_f(z)$
означает комплексную дилатацию функции $f(z),$ а $\mu_g(z)$ --
комплексную дилатацию $g.$ Согласно \cite[(1), п.~C, гл.~I]{A} для
почти всех $z\in D$ получаем:
$$f_z=\varphi_z(g(z))g_z,\qquad f_{\overline{z}}=\varphi_z(g(z))g_{\overline{z}},$$
$$\mu_f(z)=\mu_g(z)=:\mu(z), \quad K_{\mu_f}(z)=K_{\mu_g}(z):=K_{\mu}(z)=\frac{1+|\mu|}{1-|\mu|}\,.$$
Таким образом,  $K_{\mu}(z)\in L_{loc}^1(D).$ Поскольку $f$ --
конечного искажения, $g$ также конечного искажения и при почти всех
$z\in D$ выполнены соотношения
$$|\partial g|\le |\partial g|+ |\overline{\partial} g|= K^{1/2}_{\mu}(z)(|J(f, z)|)^{1/2}\,,$$
откуда по неравенству Гёльдера $|\partial g|\in L_{loc}^1 (D)$ и
$|\overline{\partial} g|\in L_{loc}^1 (D).$ Следовательно, $g\in
W_{loc}^{1,1}(D)$ и $g$ имеет конечное искажение.
\end{proof}$\Box$

\medskip
Для фиксированных области $D\subset {\Bbb C},$ числа $c\in {\Bbb C}$
и измеримой по Лебегу функции $Q:D\rightarrow [0, \infty]$ обозначим
символом $\frak{H}_{c, Q}(D)$ семейство всех гомеоморфизмов
$f:D\rightarrow {\Bbb C}\setminus\{c\}$ класса $W_{loc}^{1,1}(D)$ и
имеющих конечное искажение, таких что $K_{\mu_f}(z)\le Q(z)$ при
почти всех $z\in D.$ Для установления основных результатов настоящей
работы приведём также следующее утверждение.

\medskip
\begin{lemma}\label{lem2}

{\sl 1. Класс $\frak{H}_{a, Q}(D)$ образует нормальное семейство
отображений, как только функция $Q$ удовлетворяет одному из
следующих условий: 1) $Q\in FMO(z_0)$ в каждой точке $z_0\in D;$ 2)
$q_{z_0}(r)\le C(z_0)\cdot \log\frac{1}{r}$ при $r\rightarrow 0$ и
каждой точке $z_0\in D,$ где $C(z_0)>0$ -- некоторая постоянная; 3)
при некотором $\varepsilon_0=\varepsilon_0(z_0)$ имеет место
соотношение (\ref{eq12}). Нормальность необходимо интерпретировать в
смысле хордальной метрики $h.$

2. Если последовательность $f_n:D\rightarrow {\Bbb C},$ $f_n\in
\frak{H}_{a, Q}(D),$ сходится локально равномерно в $D$ к
отображению $f:D\rightarrow \overline{{\Bbb C}}$ в смысле метрики
$h,$ и, кроме того, функция $Q$ удовлетворяет хотя бы одному из
указанных выше условий 1)--3), то $f$ либо -- гомеоморфизм
$f:D\rightarrow {\Bbb C},$ либо -- постоянная $f:D\rightarrow
\overline{\Bbb C}.$}
\end{lemma}

\medskip
\begin{proof} Пусть $f_m\in \frak{H}_{a, Q}(D)$ -- произвольная
последовательность, тогда комплексная дилатация $\mu_{f_m}(z)$
отображения $f_m\in \frak{H}_{a, Q}(D)$ удовлетворяет следующему
условию: $|\mu_{f_m}(z)|\ne 1$ п.в., поскольку по условию леммы
$f_m$ -- конечного искажения. Поскольку $f_m$ -- гомеоморфизмы, то
либо $|\mu_{f_m}(z)|> 1$ п.в. (что соответствует случаю $J(f_m,
z)<0$ п.в.), либо $|\mu_{f_m}(z)|<1$ п.в. (что соответствует случаю
$J(f_m, z)>0$ п.в.); см., напр., \cite[разд.~V.2.2,
соотношение~(68), с.~332]{RR} либо \cite[лемма~2.14 и комментарии
после леммы~2.11]{MRV$_1$}. Не ограничивая общности рассуждений,
можно считать, что $|\mu_{f_m}(z)|< 1$ почти всюду (например,
рассмотрев для отображений $f_m,$ для которых $|\mu_{f_m}(z)|>1,$
вспомогательное семейство $\psi_m=\psi\circ f_m,$ где
$\psi(z)=x-iy,$ $z=x+iy$). В таком случае, 1-я часть заключения
леммы \ref{lem2} есть прямое следствие результатов работы
\cite[теоремы~5.1, 5.2 и следствие~5.3]{LSS}.

\medskip
Вторая часть утверждения леммы вытекает на основании
\cite[лемма~3.1]{LSS} и \cite[теоремы~4.1--4.2]{RSS} (см. также
\cite[теорема~5.3, разд.~5.3, гл.~5 и теоремы~1.3, 1.4,
следствие~1.8, разд.~1.5, гл.~1]{KSS}). $\Box$

\medskip
{\bf 4. Доказательство основных результатов.} {\it Доказательство
теоремы \ref{th1A}} основано на так называемом представлении
Стоилова. Пусть $f_m\in\frak{G}_{a, b, Q}(D)$ -- произвольная
последовательность отображений семейства $\frak{G}_{a, b, Q}(D),$
тогда согласно представлению Стоилова \cite[п.~5 (III), гл.~V]{St}
каждое отображение $f_m\in\frak{G}_{a, b, Q}(D)$ имеет вид
$f_m=\varphi_m\circ g_m,$ где $g_m$ -- некоторый гомеоморфизм, а
$\varphi_m$ -- аналитическая функция.

\medskip
Пусть $z_1, z_2$ -- две произвольные различные точки области $D.$
Рассмотрим отображения $$\Psi_m(z)=z-g_m(z_1),\quad
\psi_m(z)=\frac{z}{|g_m(z_2)-g_m(z_1)|}e^{-i\arg(g_m(z_2)-g_m(z_1))
}\,,$$
тогда $$f_m(z)=\varphi_m\circ \Psi_m^{\,-1}\circ \Psi_m\circ
g_m(z)=\varphi_m\circ
\Psi_m^{\,-1}\circ\psi_m^{\,-1}\circ\psi_m\circ\Psi_m\circ
g_m(z)\,.$$
Обозначая через $$A_m(w):=\varphi_m\circ
\Psi_m^{\,-1}\circ\psi_m^{\,-1}(w)$$ и
$$B_m(z)=\psi_m\circ\Psi_m\circ
g_m(z)=\frac{g_m(z)-g_m(z_1)}{|g_m(z_2)-g_m(z_1)|}e^{-i\arg(g_m(z_2)-g_m(z_1)}\,,$$
мы видим, что $f_m(z)=A_m\circ B_m(z),$ где $A_m$ -- аналитические
функции и $B_m$ -- гомеоморфизмы, удовлетворяющие условиям
$B_m(z_1)=0,$ $B_m(z_2)=1.$ Стоит отметить, что ввиду леммы
\ref{lem1} семейство гомеоморфизмов $\{B_m\}_{m=1}^{\infty}$
принадлежит классу $\frak{H}_{0, Q}(D\setminus\{z_1\})$ (см.
обозначения леммы \ref{lem2}), а также классу $\frak{H}_{1,
Q}(D\setminus\{z_2\}).$ Тогда ввиду первой части леммы \ref{lem2}
семейство отображений $B_m$ является нормальным семейством
отображений как в $D\setminus\{z_1\},$ так и в $D\setminus\{z_2\}.$
Поскольку произвольный компакт $C\subset D$ может быть представлен в
виде объединения $C=C_1\cup C_2,$ где $C_1$ -- компакт в
$D\setminus\{z_1\}$ и $C_2$ -- компакт в $D\setminus\{z_1\},$ то
отсюда следует, что $B_m$ также образует нормальное семейство
отображений в области $D.$

\medskip
Итак, пусть $h(B_{m_k}(x), B(x))\rightarrow 0$ при $m\rightarrow
\infty$ локально равномерно в $D,$ где $B$ -- некоторое непрерывное
отображение. Тогда, в силу условий $B_{m_k}(z_1)=0,$
$B_{m_k}(z_2)=1$ и ввиду второй части леммы \ref{lem2} отображение
$B$ является гомеоморфизмом из $D$ в ${\Bbb C}.$

\medskip
Заметим, что $B(D)\subset B_{m_k}(D)$ при всех $k\ge K_0$ и
некотором $K_0\in {\Bbb N}$ (см. \cite[предложение~1.5,
гл.~1]{KSS}). В таком случае, все отображения $A_{m_k}$ определены в
области $B(D).$ Отметим, что в этой области каждое отображение
$A_{m_k}$ не может принимать ни значение $a,$ ни значение $b,$
поскольку, в противном случае, и сами отображения $f_{m_k}$
принимали бы всё те же значения, что противоречит условию теоремы. В
таком случае, семейство отображений $A_{m_k}$ является нормальным
ввиду теоремы Монтеля о нормальности семейств аналитических функций,
не принимающих пару комплексных значений (см. \cite[$\S$ 32, гл.
II]{Mont}).

\medskip
Пусть $A_{m_{k_l}}$ -- последовательность аналитических функций,
являющаяся подпоследовательностью последовательности $A_{m_k},$
сходящаяся локально равномерно в $B(D)$ при $l\rightarrow\infty$ к
аналитической функции (либо -- тождественной бесконечности), которую
мы обозначим через $A(x).$ Пусть $C$ -- произвольный компакт в
области $D,$ тогда ввиду локально равномерной сходимости $B_{m_k}$ к
отображению $B$ все точки $B_{m_{k_l}}(x)$ лежат внутри некоторого
компакта $C_1\subset B(D)$ при всех $x\in C$ и всех $l\ge K_1\in
{\Bbb N}.$ Тогда при тех же $x$ и $l$ будем иметь:
$$h(A_{m_{k_l}}\circ B_{m_{k_l}}(x), A\circ B(x))\le$$
$$h(A_{m_{k_l}}\circ B_{m_{k_l}}(x), A\circ B_{m_{k_l}}(x))+
h(A\circ B_{m_{k_l}}(x), A\circ B(x))\le$$
$$\sup\limits_{y\in C_1}h(A_{m_{k_l}}(y), A(y))+
h(A\circ B_{m_{k_l}}(x), A\circ B(x))\rightarrow 0$$
при $l\rightarrow \infty$ равномерно по $x\in C.$ Таким образом,
последовательность $f_m\in\frak{G}_{a, b, Q}(D)$ имеет
подпоследовательность, сходящуюся локально равномерно в $D.$
\end{proof}$\Box$

\medskip
{\it Доказательство следствия \ref{cor1A}} вытекает из того, что
семейство $\frak{G^{\,*}}_{a, b, B}(D)$ является подклассом
семейства $\frak{G}_{a, b, B}(D)$ при указанных условиях на функцию
$B.$ Действительно, $\frak{G^{\,*}}_{a, b, B}(D)\subset W_{loc}^{1,
1}$ и каждое $f\in\frak{G^{\,*}}_{a, b, B}(D)$ имеет конечное
искажение ввиду \cite[следствия~3.3 и 3.5]{SS}, кроме того,
$K_{\mu}(z)\le Q(z)$ почти всюду ввиду \cite[следствие~3.2]{SS}.
$\Box$

{\bf 5. О компактности классов Соболева.} Для фиксированных области
$D\subset {\Bbb C},$ чисел $a, b\in D,$ $a\ne b,$ $a^{\,\prime},
b^{\,\prime}\in {\Bbb C},$ $a^{\,\prime}\ne b^{\,\prime},$ и
измеримой по Лебегу функции $Q:D\rightarrow [0, \infty],$ обозначим
символом $\frak{A}_{a, b, a^{\,\prime}, b^{\,\prime}, Q}(D)$
семейство всех открытых дискретных отображений $f:D\rightarrow {\Bbb
C}$ класса $W_{loc}^{1,1}(D)$ и имеющих конечное искажение, таких
что $f(a)=a^{\,\prime},$ $f(b)=b^{\,\prime}$ и $K_{\mu_f}(z)\le
Q(z)$ при почти всех $z\in D.$ Справедливо следующее утверждение.

\medskip
\begin{theorem}\label{th1}
{\sl Класс $\frak{A}_{a, b, a^{\,\prime}, b^{\,\prime}, Q}(D)$
является компактным (т.е., нормальным и замкнутым семейством
отображений в топологии локально равномерной сходимости)}.
\end{theorem}

\medskip
\begin{proof} Заметим, что семейство отображений является нормальным
ввиду теоремы \ref{th1A}. Более того, повторяя доказательство этой
теоремы, мы приходим к заключению, что произвольная сходящаяся
последовательность $f_m$ представима в виде композиции $f_m=A_m\circ
B_m,$ где $B_m$ -- последовательность гомеоморфизмов, сходящаяся
локально равномерно к гомеоморфизму $B,$ а $A_m$ --
последовательность аналитических функций, сходящаяся к аналитической
функции $A.$ При этом, $f_m\rightarrow f:=A\circ B.$ Осталось
показать, что $f\in \frak{A}_{a, b, a^{\,\prime}, b^{\,\prime},
Q}(D).$

\medskip
Заметим, что условия нормировки $f_m(a)=a^{\,\prime},$
$f_m(b)=b^{\,\prime}$ исключают возможность, когда $A$ является
постоянной функцией. Значит, $f$ дискретно и открыто. Кроме того,
ввиду леммы \ref{lem1} заключаем, что $B_m\in W_{loc}^{1, 1},$ $B_m$
имеют конечное искажение и $K_{\mu_{B_m}}(z)=K_{{\mu}_{f_m}}(z).$
Ввиду теорем \cite[17.1--17.2]{RSS} предельное отображение $B$ также
принадлежит классу $W_{loc}^{1,1},$ имеет конечное искажение и его
максимальная дилатация $K_{\mu_B}(z)$ не превосходит $Q(z)$ почти
всюду. Тогда, очевидно, $f=A\circ B\in W_{loc}^{1,1}$ и
$K_{\mu_f}(z)\ne Q(z).$ Равенства $f(a)=a^{\,\prime}$ и
$f(b)=b^{\,\prime}$ элементарно получаются предельным переходом по
$m$ из равенств $f_m(a)=a^{\,\prime}$ и $f_m(b)=b^{\,\prime}.$
Теорема доказана. $\Box$
\end{proof}

\noindent{{\bf Евгений Александрович Севостьянов} \\ Институт
прикладной математики и механики НАН Украины \\
83 114 Украина, г. Донецк, ул. Розы Люксембург, д. 74, \\
тел. +38 (066) 959 50 34 (моб.), +38 (062) 311 01 45 (раб.), e-mail:
brusin2006@rambler.ru, esevostyanov2009@mail.ru}

\end{document}